\documentclass[letterpaper, 11pt,  reqno]{amsart}

\usepackage{amsmath,amssymb,amscd,amsthm,amsxtra, esint}

\usepackage[implicit=true]{hyperref}

\allowdisplaybreaks[2]

\usepackage{color}

\definecolor{gr}{rgb}   {0.,   0.69,   0.23 }
\definecolor{bl}{rgb}   {0.,   0.5,   1. }
\definecolor{mg}{rgb}   {0.85,  0.,    0.85}
\definecolor{yl}{rgb}   {0.8,  0.7,   0.}
\definecolor{or}{rgb}  {0.7,0.2,0.2}

\renewcommand{\H}{\mathcal{H}}

\newtheorem{theorem}{Theorem} [section]

\newtheorem{lemma}[theorem]{Lemma}

\newtheorem{remark}[theorem]{Remark}

\newtheorem*{acknowledgment}{Acknowledgments}

\DeclareMathOperator*{\supp}{supp}

\newcommand{\1}{\hspace{0.5mm}\text{I}\hspace{0.5mm}}
\newcommand{\II}{\text{I \hspace{-2.8mm} I} }
\newcommand{\III}{\text{I \hspace{-2.9mm} I \hspace{-2.9mm} I}}

\newcommand{\noi}{\noindent}
\newcommand{\Z}{\mathbb{Z}}
\newcommand{\R}{\mathbb{R}}
\newcommand{\C}{\mathbb{C}}
\newcommand{\T}{\mathbb{T}}

\let\P= \undefined
\newcommand{\P}{\mathbf{P}}

\newcommand{\M}{\mathcal{M}}

\newcommand{\N}{\mathcal{N}}
\newcommand{\NB}{\mathbb{N}}

\newcommand{\F}{\mathcal{F}}

\newcommand{\al}{\alpha}

\newcommand{\nb}{\nabla}

\newcommand{\Dl}{\Delta}
\newcommand{\eps}{\varepsilon}

\newcommand{\g}{\gamma}
\newcommand{\G}{\Gamma}
\newcommand{\ld}{\lambda}
\newcommand{\Ld}{\Lambda}
\newcommand{\s}{\sigma}
\newcommand{\Si}{\Sigma}
\newcommand{\ft}{\widehat}

\newcommand{\wt}{\widetilde}
\newcommand{\cj}{\overline}
\newcommand{\dx}{\partial_x}
\newcommand{\dt}{\partial_t}

\renewcommand{\l}{\ell}
\renewcommand{\o}{\omega}
\renewcommand{\O}{\Omega}

\newcommand{\les}{\lesssim}
\newcommand{\ges}{\gtrsim}

\newcommand{\jb}[1]
{\langle #1 \rangle}

\newcommand{\ind}{\mathbf 1}

\numberwithin{equation}{section}
\numberwithin{theorem}{section}

\usepackage{tikz} 
\usepackage{marginnote}
\usepackage{scalerel} % to rescale the paraproduct symbols

\usetikzlibrary{shapes.misc}
\usetikzlibrary{shapes.symbols}
\usetikzlibrary{decorations}
\usetikzlibrary{decorations.markings}

\tikzset{
	dot/.style={circle,fill=black,draw=black,inner sep=0pt,minimum size=0.5mm},
	>=stealth,
	}

\tikzset{
	dot2/.style={circle,fill=black,draw=black,inner sep=0pt,minimum size=0.2mm},
	>=stealth,
	}

\tikzset{
	ddot/.style={circle,fill=white,draw=black,inner sep=0pt,minimum size=0.8mm},
	>=stealth,
	}

\tikzset{decision/.style={ % requires library shapes.geometric
        draw,
        diamond,
        aspect=1.5
    }}

\tikzset{dia2/.style
={diamond,fill=white,draw=black,inner sep=0pt,minimum size=1mm},
	>=stealth,
	}

\tikzset{dia/.style
={star,fill=black,draw=black,inner sep=0pt,minimum size=1mm},
	>=stealth,
	}

\tikzset{dia/.style
={diamond,fill=black,draw=black,inner sep=0pt,minimum size=1.3mm},
	>=stealth,
	}

\makeatletter
\def\DeclareSymbol#1#2#3{\expandafter\gdef\csname MH@symb@#1\endcsname{\tikz[baseline=#2,scale=0.15]{#3}}}
\def\<#1>{\csname MH@symb@#1\endcsname}
\makeatother
\newcommand{\pe}{\mathbin{\scaleobj{0.7}{\tikz \draw (0,0) node[shape=circle,draw,inner sep=0pt,minimum size=8.5pt] {\scriptsize  $=$};}}}

\newcommand{\pl}{\mathbin{\scaleobj{0.7}{\tikz \draw (0,0) node[shape=circle,draw,inner sep=0pt,minimum size=8.5pt] {\scriptsize $<$};}}}
\newcommand{\pg}{\mathbin{\scaleobj{0.7}{\tikz \draw (0,0) node[shape=circle,draw,inner sep=0pt,minimum size=8.5pt] {\scriptsize $>$};}}}

\newcommand*\circled[1]{\tikz[baseline=(char.base)]{
            \node[shape=circle,draw,inner sep=0.5pt] (char) {#1};}}

\begin{document}
\baselineskip = 14pt

\title[On the probabilistic Cauchy theory 
for nonlinear dispersive PDEs]
{On the probabilistic Cauchy theory 
for nonlinear dispersive PDEs}

\author[\'A.~B\'enyi, T.~Oh, and O.~Pocovnicu]
{\'Arp\'ad  B\'enyi, Tadahiro Oh, and Oana Pocovnicu}

\address{
\'Arp\'ad  B\'enyi\\
Department of Mathematics\\
 Western Washington University\\
 516 High Street, Bellingham\\
  WA 98225\\ USA}
\email{arpad.benyi@wwu.edu}

\address{
Tadahiro Oh\\
School of Mathematics\\
The University of Edinburgh\\
and The Maxwell Institute for the Mathematical Sciences\\
James Clerk Maxwell Building\\
The King's Buildings\\
 Peter Guthrie Tait Road\\
Edinburgh\\ 
EH9 3FD\\United Kingdom} 

\email{hiro.oh@ed.ac.uk}

\address{
Oana Pocovnicu\\
Department of Mathematics, Heriot-Watt University and The Maxwell Institute for the Mathematical Sciences, Edinburgh, EH14 4AS, United Kingdom}
\email{o.pocovnicu@hw.ac.uk}



\begin{abstract}
In this note, we review some of the recent 
developments in the well-posedness theory of nonlinear dispersive partial differential equations %(PDEs)
with random initial data.

\end{abstract}

%\date{\today}

%%
%
\maketitle
%
%\tableofcontents

\baselineskip = 13pt

\section{Introduction}
Nonlinear dispersive partial differential equations (PDEs) naturally appear as models describing wave phenomena in various branches of physics and engineering such as quantum mechanics, nonlinear optics, plasma physics, water waves, and atmospheric sciences. They have received wide attention from the applied 
science community due to their importance in applications and have also been studied extensively from the theoretical point of view, providing a framework for development of analytical ideas and tools. 

The simplest, yet most important examples of nonlinear dispersive PDEs
are the following
nonlinear Schr\"odinger equations (NLS):\footnote{For conciseness, 
we restrict our attention to the defocusing case
in the following.}
\begin{equation}
\begin{cases}\label{NLS1}
i \partial_t u =  \Delta u - |u|^{p-1} u  \\
u|_{t = 0} = u_0, 
\end{cases}
\qquad ( t, x) \in \R \times \M
\end{equation}

\noi
and 
nonlinear wave equations (NLW):
\begin{equation}
\begin{cases}\label{NLW1}
 \dt^2  u = \Delta u - |u|^{p-1}u  \\
(u, \dt u) |_{t = 0} = (u_0, u_1), 
\end{cases}
\qquad ( t, x) \in \R \times \M,
\end{equation}

\noi
where $\M = \R^d$ or $\T^d$  and $p > 1$.
Over the last several decades, 
multilinear harmonic analysis has played a crucial role in building basic insights on the study of nonlinear dispersive PDEs, settling questions on existence of solutions to these equations, their long-time behavior, and singularity formation. 
Furthermore, in recent years, a remarkable combination of PDE techniques and probability theory has had a significant impact on the field. 
In this note, we go over some of the recent developments in this direction.

In the classical deterministic well-posedness theory, 
the main goal is to construct unique solutions for {\it all} initial data belonging to a certain 
fixed function space
such as the $L^2$-based Sobolev spaces: 
\[ \text{$H^s(\M)$ for \eqref{NLS1}
\quad and \quad 
$\H^s(\M) : = H^s(\M) \times H^{s-1}(\M)$ for \eqref{NLW1}.}\]

\noi 
In practice, however, we are often interested in the typical behavior of solutions. 
Namely, even if certain pathological behavior occurs, we may be content if we can show that almost all solutions behave well and do not exhibit such pathological behavior. 
This concept may be formalized in terms of probability. 
For example, 
in terms of well-posedness theory, 
one may consider an evolution equation with random initial data and try to construct unique solutions in an almost sure manner.
This idea first appeared in 
Bourgain's seminal paper
\cite{BO94},
where he constructed 
%
%the probabilistic construction of 
%
global-in-time solutions
to NLS on $\T$ 
almost surely with respect to the random initial data distributed according to the Gibbs measure.
See Subsection \ref{SUBSEC:GWP}.

Such probabilistic construction of solutions with random initial data also allows
us to go beyond deterministic thresholds in certain situations. 
First recall that 
NLS \eqref{NLS1} 
and NLW \eqref{NLW1}
 enjoy the following dilation symmetry:
\begin{align}
%\begin{split}
 u(t, x) \longmapsto u^\ld(t, x) = \ld^{-\frac{2}{p-1}} u (\ld^{-\al}t, \ld^{-1}x),
% u(t, x) \longmapsto u^\ld(t, x) = \ld^{-\frac{2}{p-1}} u (\ld^{-2}t, \ld^{-1}x).
%\end{split}
 \label{scaling1}
\end{align}

\noi
with $\al = 2$ for \eqref{NLS1} and $\al = 1$ for \eqref{NLW1}.
Namely, if $u$ is a solution to \eqref{NLS1} or \eqref{NLW1} 
on $\R^d$, then $u^\ld$ is also a solution 
to the same equation on $\R^d$ with the rescaled initial data.
This dilation symmetry induces the following  scaling-critical Sobolev regularity:
 \begin{align}
 s_\text{crit} = \frac d2 - \frac{ 2}{p-1}
\label{scaling2}
 \end{align}
 
 \noi
such that the homogeneous $\dot{H}^{s_\text{crit}}(\R^d)$-norm is invariant
under the dilation symmetry.
This critical regularity $s_\text{crit}$ provides
a threshold regularity for well-posedness and ill-posedness
of \eqref{NLS1} and \eqref{NLW1}.\footnote{In fact, there are other
 critical regularities induced by the Galilean invariance for \eqref{NLS1}
 and the Lorentzian symmetry for \eqref{NLW1}
 below which the equations are ill-posed; see \cite{LS, CCT2, MOLI,  GuoOh}.
 We point out, however, that these additional critical regularities
 are relevant only when the dimension is low and/or the degree $p$ is small. 
For example, for NLS \eqref{NLS1} with an algebraic nonlinearity ($p \in 2\NB + 1$), 
the critical regularity induced by the Galilean invariance is  relevant (i.e.~higher than 
the scaling-critical regularity $s_\text{crit}$ in \eqref{scaling2}) only for $d = 1$
and $p = 3$.
 For simplicity, we only consider the scaling-critical regularities in the following.
 }
While there is no dilation symmetry in the periodic setting, 
the heuristics provided by the scaling argument also plays an important role.
On the one hand, 
there is a good local-in-time theory 
for \eqref{NLS1} and \eqref{NLW1}
(at least when the dimension $d$ and the degree $p$
are not too small).
See \cite{BOP2, POC, O17} for the references therein.
On the other hand, 
it is known that \eqref{NLS1} and \eqref{NLW1} are ill-posed
in the supercritical regime:  $s < s_\text{crit}$.
See \cite{CCT, BT1, Kishimoto, 
OW, Xia, O17, CP, OOTz}. 
Regardless of these ill-posedness results, 
by considering random initial data (see Section \ref{SEC:data}), 
one may still prove almost sure local well-posedness\footnote{Namely, local-in-time
existence of unique solutions almost surely with respect to given random initial data.}
in the supercritical regime. %at  the supercritical regularities  ($s < s_\text{crit}$).
This probabilistic construction of local-in-time solutions was first implemented
by Bourgain \cite{BO96} 
in the context of NLS
and by Burq-Tzvetkov \cite{BT1}
in the context of NLW.
In more recent years, there have also been 
results on  almost sure global well-posedness
for these equations; 
see for example
\cite{CO, BT3, LM, POC, OP, SX, KMV, OOP}.
See also the lecture note by Tzvetkov \cite{TzNote}.
We will discuss some  aspects of probabilistic well-posedness  in Section \ref{SEC:WP}.

\section{On random initial data}\label{SEC:data}
In this section, 
we go over random initial data 
based on 
random Fourier series on $\T^d$
and its analogue on $\R^d$.

\subsection{Random initial data on $\T^d$}
In the context of nonlinear dispersive PDEs, 
probabilistic construction of solutions was initiated
in an effort to construct
well-defined dynamics
almost surely with respect to the Gibbs measure
for NLS on $\T^d$, $d = 1, 2$
\cite{BO94, McKean, BO96}.
Before discussing this problem for NLS on $\T^d$, 
let us consider the following 
finite dimensional  Hamiltonian flow on $\mathbb{R}^{2N}$:
\begin{equation} \label{HR2}
\dot{p}_n = \frac{\partial H}{\partial q_n}
\qquad \text{and}\qquad  
\dot{q}_n = - \frac{\partial H}{\partial p_n},  
\end{equation}

\noi
$n = 1, \dots, N$, 
 with Hamiltonian $ H (p, q)= H(p_1, \cdots, p_N, q_1, \cdots, q_N)$.
Liouville's theorem states that the Lebesgue measure 
$d p dq = \prod_{n = 1}^N dp_n dq_n$
on $\mathbb{R}^{2N}$ is invariant under the flow.
Then, it follows from the conservation of the Hamiltonian $H(p, q)$
that  the Gibbs measure:\footnote{Hereafter, we use $Z$ to denote various normalizing constants
so that the resulting measure is a probability measure provided that it makes sense.}
\[d\rho = Z^{-1} e^{- H(p, q)} dp dq\]

\noi
is  invariant under the flow of \eqref{HR2}.
Recall that NLS \eqref{NLS1}  is a Hamiltonian PDE
with the following Hamiltonian:
\begin{align}
H(u) = \frac{1}{2} \int_\M |\nb  u |^2 dx + \frac{1}{p+1}\int_\M |u|^{p+1} dx. 
\label{Hamil1}
\end{align}

\noi
Moreover, the mass $M(u)$ defined by 
\begin{align}
 M(u) = \frac{1}{2}\int_\M|u|^2 dx
\label{mass1}
 \end{align}

\noi
is conserved under the dynamics of \eqref{NLS1}.
Then, by drawing an analogy to the finite dimensional case, 
one may expect that the Gibbs measure:\footnote{Here, we added the mass in the exponent to avoid
a problem at the zeroth frequency in \eqref{gauss2} below.}
\begin{align}
\text{``}d\rho = Z^{-1} e^{- H(u)- M(u)} du\text{''}
\label{Gibbs1}
\end{align}

\noi
is invariant under the flow of \eqref{NLS1}.
Here, the expression in \eqref{Gibbs1} is merely formal,
where ``$du$'' denotes the non-existent Lebesgue measure
on an infinite dimensional phase space.

We first introduce a family of mean-zero Gaussian measures $\mu_s$, $s \in \R$,  on 
periodic distributions on $\T^d$, formally given by
\begin{align}
 d \mu_s 
  = Z^{-1} e^{-\frac 12 \| u\|_{H^s}^2} du 
  = Z^{-1} \prod_{n \in \Z^d} e^{-\frac 12 \jb{n}^{2s} |\ft u(n)|^2}   d\ft u(n), 
\label{mu}
 \end{align}

\noi
where  $\jb{\,\cdot\,} = (1+|\cdot|^2)^\frac{1}{2}$.
As we see in \eqref{Gibbs2} below, the Gibbs measure $\rho$ is constructed
as the Gaussian measure $\mu_1$ with a weight.
While the expression  
$d \mu_s 
  = Z^{-1} \exp(-\frac 12 \| u\|_{H^s}^2 )du $
may suggest that $\mu_s$ is a Gaussian measure on $H^s(\T^d)$, 
we need to enlarge a space in order to make sense of $\mu_s$
as a countably additive probability measure.
In fact, the Gaussian measure $\mu_s$ defined above is  the induced probability measure
under the following map:\footnote{In the following, we drop the harmless factor of $2\pi$.}
\begin{align}
 \o \in \O \longmapsto u^\o(x) = u(x; \o) = \sum_{n \in \Z^d} \frac{g_n(\o)}{\jb{n}^s}e^{in\cdot x}, 
\label{gauss1}
 \end{align}

\noi
where
$\{ g_n \}_{n \in \Z^d}$ is a sequence of independent standard complex-valued 
Gaussian random variables on a probability space $(\O, \F, P)$. % i.e.~$\text{Var}(g_n) = 2$.
It is easy to check that the random function $u^\o$ in \eqref{gauss1}
 lies  in 
 \[ \text{$H^{\s}(\T^d)$ for $\s < s -\frac d2$
\quad but \quad not in $H^{s-\frac d2}(\T^d)$},\] almost surely. 
Moreover,  for the same range of $\s$,  
 $\mu_s$ is a Gaussian probability measure on $H^{\s}(\T^d)$ 
and the triplet $(H^s(\T^d), H^\s(\T^d), \mu_s)$ forms an abstract Wiener space.
See \cite{GROSS, Kuo2}.
Note that, when $s = 1$, the random Fourier series \eqref{gauss1}
basically corresponds to the
Fourier-Wiener series for the Brownian motion.
See \cite{BenyiOh} for more on this.

By restricting our attention to $\T^d$, 
we substitute the expressions \eqref{Hamil1} and \eqref{mass1} in \eqref{Gibbs1}.
Then, 
we formally obtain 
\begin{align}
d\rho 
& = Z^{-1} e^{- \frac{1}{p+1} \int |u|^{p+1} } e^{-\frac{1}{2} \| u\|_{H^1}^2} du\notag\\
& = Z^{-1} e^{- \frac{1}{p+1} \int |u|^{p+1} } d\mu_1, 
\label{Gibbs2}
\end{align}

\noi
where $\mu_1$ is as in \eqref{mu} with $s = 1$.
When $d = 1$, it is easy to see that the Gibbs measure $\rho$
is a well-defined probability measure, 
absolutely continuous with respect to the Gaussian measure $\mu_1$.
In particular, it is supported on $H^\s(\T)$ for any $\s < \frac 12$.
When $d = 2$, a typical element under $\mu_1$ lies in $H^{-\eps}(\T^2) \setminus L^2(\T^2)$
for any $\eps > 0$.
As such, the weight
$e^{- \frac{1}{p+1} \int |u|^{p+1} }$ in \eqref{Gibbs2} equals 0 almost surely
and  hence the expression \eqref{Gibbs2} for $\rho$
does not make sense as a probability measure.
Nonetheless, when $p \in 2\NB + 1$, one can apply a suitable renormalization
(the Wick ordering) and construct the Gibbs measure $\rho$
associated with the Wick ordered Hamiltonian
such that $\rho$ is absolutely continuous to $\mu_1$.\footnote{When $ d\geq 3$, it is known that the Gibbs measure $\rho$ can be constructed
only for $d = 3$ and $p = 3$.
In this case, the resulting Gibbs measure $\rho$ is not absolutely continuous with respect to the Gaussian measure $\mu_1$.  See \cite{AK} for the references therein, regarding the construction of the Gibbs measure (the $\Phi^4_3$ measure) in the real-valued setting.}
See \cite{BO96, OTh, OTh2}
for more on the renormalization in the two-dimensional case.
This shows that when $d = 1, 2$, 
it is of importance to study the dynamical property of NLS \eqref{NLS1}
with the random initial data $u_0^\o$ distributed according to the Gaussian measure $\mu_1$,
namely given by the random Fourier series \eqref{gauss1} with $s = 1$:
\begin{align}
u_0^\o(x) =  \sum_{n \in \Z^d} \frac{g_n(\o)}{\jb{n}}e^{in\cdot x}
\in H^{1 - \frac d 2 - \eps}(\T^d) \setminus H^{1 - \frac d2}(\T^d)
\label{gauss2}
 \end{align}

\noi
almost surely for any $\eps > 0$.
In \cite{BO96}, 
Bourgain studied the (renormalized) cubic NLS on $\T^2$ with the random initial data \eqref{gauss2}.
Recalling that the two-dimensional cubic NLS is $L^2$-critical, 
we see that this random initial data $u_0^\o$ lies slightly below the critical regularity.\footnote{In terms
of the Besov spaces  $B^\s_{p, \infty}$, $p < \infty$, 
we see that $u_0^\o$ in \eqref{gauss2} lies almost surely in the critical spaces 
 $B^0_{2, \infty}$.  See \cite{BenyiOh}.}
Nonetheless, by combining the deterministic analysis 
(the Fourier restriction norm method introduced in \cite{BO1})
with the probabilistic tools, in particular, the probabilistic Strichartz estimates (Lemma \ref{LEM:PStr1}), 
he managed to prove almost sure local well-posedness
with respect to the random initial data $u_0^\o$ in \eqref{gauss2}.

In the context of the cubic NLW on a three-dimensional compact Riemannian manifold $\M$, 
Burq-Tzvetkov \cite{BT1} considered the Cauchy problem with a more general class of random initial data.
In particular, given a rough initial data $(u_0, u_1) \in \H^s(\M)$
with $s < s_\text{crit} = \frac 12$,
they introduced a randomization  $(u_0^\o, u_1^\o)$
of the given initial data $(u_0, u_1)$
and proved almost sure local well-posedness
with respect to this randomization. 
For simplicity, we discuss this randomization for a single function $u_0$ on $\T^d$ in the following.
Fix $u_0 \in H^s(\T^d)$.
Then, we define a randomization 
$u_0^\o$ of $u_0$ by 
setting
\begin{align}
u_0^\o(x) :=  \sum_{n \in \Z^d} g_{ n}(\o)\ft u_0(n)e^{in \cdot x}, 
\label{gauss3}
 \end{align}

\noi
where $\ft u_0(n)$ denotes the Fourier coefficient of $u_0$
and 
$\{ g_{n} \}_{n \in \Z^d}$ is a sequence of independent mean-zero 
complex-valued random variables
with bounded moments up to a certain order.\footnote{In the real-valued setting, 
we also need to impose that $g_{ -n} = \cj{g_{n}}$ so that,
given a real-valued function $u_0$,  the resulting randomization
$u_0^\o$ remains real-valued.
A similar comment applies to the randomization \eqref{gauss4}
introduced for functions on $\R^d$.}
Note that the random Fourier series in \eqref{gauss1} and \eqref{gauss2}
can be viewed as a randomization of the particular function $u_0$ with the Fourier coefficient $\jb{n}^{-s}$
by independent standard Gaussian random variables $\{ g_n\}_{n \in \Z^d}$.
The main point of the randomization \eqref{gauss3}
is that while the randomized function $u_0^\o$ does not
enjoy any smoothing in terms of differentiability, 
it enjoys a gain of integrability (Lemma \ref{LEM:Hs}).

\subsection{Probabilistic Strichartz estimates}
\label{SUBSEC:PStr}

In this subsection, 
we discuss the effect of the randomization \eqref{gauss3} introduced in the previous subsection.
For simplicity, 
we further assume that the 
 probability distributions $\mu_n^{(1)}$ and $\mu_n^{(2)}$
of  the real and imaginary parts of $g_n$ in \eqref{gauss3}
are independent
and satisfy 
\begin{equation}
\bigg| \int_{\R} e^{\g x } d \mu_n^{(j)}(x) \bigg| \leq e^{c\g^2}
\label{R2}
\end{equation}
	
\noi
for all $\g \in \R$, $n \in \Z^d$, $j = 1, 2$.
Note that \eqref{R2} is satisfied by
standard complex-valued Gaussian random variables,
Bernoulli random variables,
and any random variables with compactly supported distributions.
Under this extra assumption \eqref{R2}, we have the following estimate.
 See \cite{BT1} for the proof.

\begin{lemma} \label{LEM:R1}
Assume \eqref{R2}. Then, there exists $C>0$ such that
\[ \bigg\| \sum_{n \in \Z^d} g_n(\omega) c_n\bigg\|_{L^p(\Omega)}
\leq C \sqrt{p} \| c_n\|_{\l^2_n(\Z^d)}\]

\noi
for all $p \geq 2$ and $\{c_n\} \in \l^2(\Z^d)$.
\end{lemma}

When $\{ g_n \}_{n \in \Z^d}$ is a sequence of independent standard 
Gaussian random variables, 
Lemma \ref{LEM:R1} follows from the Wiener chaos estimate 
(see Lemma \ref{LEM:hyp} below) with $k = 1$.

Given $u_0 \in H^s(\T^d)$, it is easy to see that its randomization $u_0^\o$
in \eqref{gauss3} lies in $ H^s(\T^d)$
almost surely.
One can also show that, under some non-degeneracy condition,  there is no smoothing upon randomization
in terms of differentiability;
see, for example, Lemma B.1 in \cite{BT1}.
The main point of the randomization \eqref{gauss3} is
its improved integrability.

\begin{lemma} \label{LEM:Hs}
Given  $u_0 \in L^2(\T^d)$, let $u_0^\o$ be its randomization defined in \eqref{gauss3},
 satisfying \eqref{R2}.
Then, given finite $p \geq 2$, there exist $C, c >0$ such that
\begin{align*}
P\Big( \| u_0^\omega \|_{ L^p} > \ld\Big)
\leq C e^{-c \ld^2  \|u_0\|_{ L^2}^{-2}}
\end{align*}

\noi
for all $\ld > 0$.
In particular, $u_0^\o$ lies in $L^p(\T^d)$ almost surely.

\end{lemma}

Such gain of integrability is well known  for random Fourier series;
see for example \cite{PZ, Kahane, AT}.
The proof of Lemma \ref{LEM:Hs} is standard and follows from Minkowski's integral inequality, Lemma \ref{LEM:R1}, 
and Chebyshev's inequality.
See \cite{BT1, CO, BOP1}.
By a similar argument, one can also establish the following probabilistic improvement
of the Strichartz estimates.

\begin{lemma}\label{LEM:PStr1}
Given $u_0$ on $\T^d$, 
let $u_0^\o$ be its randomization defined in \eqref{gauss3}, satisfying \eqref{R2}.
Then,
given finite  $q \geq 2$ and $2 \leq  r \leq \infty$,
there exist $C, c>0$ such that
\begin{align*}
P\Big( \| e^{-it \Dl} u_0^\omega\|_{L^q_t L^r_x([0, T]\times \T^d)}> \ld\Big)
\leq C\exp \bigg(-c \frac{\ld^2}{ T^\frac{2}{q}\|u_0\|_{H^s}^{2}}\bigg)
\end{align*}
	
\noi
for all  $ T > 0$ and $\lambda>0$
with  \textup{(i)} $s = 0$ if $r < \infty$
and  \textup{(ii)}  $s > 0$ if $r = \infty$.

\end{lemma}

By setting $\ld = T^\theta \|u_0\|_{L^2}$, we have
\begin{equation}
\|e^{-it\Dl} u_0^\o\|_{L^q_tL^r_x([0, T]\times \T^d)}
\les T^\theta \|u_0\|_{L^2(\T^d)}
\label{PStr2}
\end{equation}

\noi
outside a set of probability at most
$ C \exp \big(-c T^{2\theta - \frac{2}{q}}\big).$
Note that, as long as $\theta < \frac{1}{q}$, this probability can be made arbitrarily small by letting $T\to 0$.
We can interpret \eqref{PStr2}
as the probabilistic improvement
of the usual Strichartz estimates, 
where the indices $q$ and $r$ 
satisfy certain relations\footnote{See  \eqref{Str2} below for the scaling condition on $\R^d$.}
and the resulting estimates come with possible loss of derivatives.
See \cite{GOW, BD}.
In Lemma \ref{LEM:PStr1}, we only stated
the probabilistic Strichartz estimates
for the Schr\"odinger equation.
Similar probabilistic Strichartz estimates also hold
for the wave equation.
See \cite{BT1, BT3, POC}.

On the one hand, the probabilistic Strichartz estimates in Lemma \ref{LEM:PStr1}
allow us to exploit the randomization at the linear level.
On the other hand, 
the following Wiener chaos estimate (\cite[Theorem~I.22]{Simon})
allows us to exploit the randomization at a multilinear level.
See also \cite[Proposition~2.4]{TTz}.

\begin{lemma}\label{LEM:hyp}
 Let $\{ g_n\}_{n \in \Z^d }$ be 
 a sequence of  independent standard 
 Gaussian random variables.
Given  $k \in \mathbb{N}$, 
let $\{P_j\}_{j \in \NB}$ be a sequence of polynomials in 
$\bar g = \{ g_n\}_{n \in \Z^d }$ of  degree at most $k$. 
Then, for $p \geq 2$, we have
\begin{equation*}
 \bigg\|\sum_{j \in \NB} P_j(\bar g) \bigg\|_{L^p(\O)} \leq (p-1)^\frac{k}{2} \bigg\|\sum_{j \in \NB} P_j(\bar g) \bigg\|_{L^2(\O)}.
 \end{equation*}

\end{lemma}

This lemma is a direct corollary to the
  hypercontractivity of the Ornstein-Uhlenbeck
semigroup due to Nelson \cite{Nelson2}.
It %Lemma \ref{LEM:hyp}
allows us to prove the following
probabilistic improvement of Young's inequality.
Such a probabilistic improvement was essential
in the probabilistic   construction of  solutions 
to the (renormalized) cubic NLS on $\T^d$, $d = 1, 2$, 
in a low regularity setting \cite{BO96, CO}.
For simplicity, we consider a trilinear case.

\begin{lemma}\label{LEM:PY}
Let $a_n, b_n, c_n \in \l^2(\Z^d; \C)$.
Given  a sequence $\{g_n \}_{n \in \Z^d}$ of independent 
standard complex-valued Gaussian random variables, 
define $a^\o_n = g_n a_n$, $b^\o_n = g_n b_n$, and $c^\o_n = g_n c_n$, $n \in \Z^d$.
Then, given $\eps > 0$, there exists a set $\O_\eps \subset \O$
with $P(\O_\eps^c) < \eps$ 
and $C_\eps> 0$ such that\footnote{One can choose $C_\eps = 	
\big(\frac{1}{c} \log \frac{C}{\eps}\big)^\frac{3}{2}$.}
 \begin{align}
\big\| a^\o_n * b_n^\o*c^\o_n\big\|_{\l^2}
\leq C_\eps 
\| a_n\|_{\l^2}\| b_n\|_{\l^2}\| c_n\|_{\l^2}
\label{B8}
\end{align}

\noi
for all $\o \in \O_\eps$.
\end{lemma}
	
The proof of Lemma \ref{LEM:PY} is immediate from the following tail estimate:
\begin{align*}
P\bigg(\Big|\sum_{j \in \NB} P_j(\bar g)\Big| >\ld \bigg) 
\leq C \exp\bigg(-c\Big\|\sum_{j \in \NB} P_j(\bar g)\Big\|_{L^2(\O)}^{-\frac{2}{k}} \ld^\frac{2}{k}\bigg),
\end{align*}

\noi
which is a consequence of Lemma \ref{LEM:hyp} and Chebyshev's inequality.
Note that Young's inequality (without randomization) only yields
 \begin{align*}
\big\| a_n * b_n*c_n \big\|_{\l^2}
\leq
\| a_n\|_{\l^2}\| b_n\|_{\l^1}\| c_n\|_{\l^1}.
\end{align*}

\noi
Recalling that $\l^1 \subset \l^2$, we see that 
there is 
 a significant improvement 
 in \eqref{B8} 
under randomization of the sequences,
which was a key in establishing crucial nonlinear estimates
in a probabilistic manner in \cite{BO96, CO}.

\subsection{Random initial data on $\R^d$}
We conclude this section by briefly going over the randomization
of a function on $\R^d$
analogous to \eqref{gauss3} on $\T^d$.  See \cite{ZF, LM, BOP1}.
On compact domains, 
 there is  a countable basis of eigenfunctions of the Laplacian
and thus there is a natural way to introduce a randomization via \eqref{gauss3}.
On the other hand, on $\R^d$,
there is no  countable basis of $L^2(\R^d)$ consisting of eigenfunctions of the Laplacian
and hence there is no ``natural'' way to introduce a randomization.
In the following, we discuss  a randomization 
adapted to the so-called  Wiener decomposition \cite{W}
of the frequency space: $\R^d = \bigcup_{n \in \Z^d} Q_n$, 
where $Q_n$ is the unit cube centered at $n \in \Z^d$.

Let $\psi \in \mathcal{S}(\R^d)$ such that
\begin{equation*}
\supp \psi \subset [-1,1]^d
\quad \text{and} \quad\sum_{n \in \Z^d} \psi(\xi -n) \equiv 1
\ \text{ for any }\xi \in \R^d.
\end{equation*}

\noi
Then, given a function $u_0$ on $\R^d$, 
we have
\begin{align*}
u_0 = \sum_{n \in \Z^d} \psi(D-n) u_0,
\end{align*}

\noi
where $\psi(D-n)$ is defined by 
 $\psi(D-n)u_0(x)=\int_{\R^d} \psi (\xi-n)\ft u_0 (\xi)e^{ ix\cdot \xi} d\xi$, 
 namely, 
 the Fourier multiplier operator with symbol $\ind_{Q_n}$ conveniently smoothed.
This decomposition leads to the following randomization of $u_0$
adapted to the Wiener decomposition.
Let $\{g_n\}_{n \in \Z^d}$ be a sequence of independent mean-zero complex-valued random variables
as in \eqref{gauss3}.
Then, we can define the  Wiener randomization\footnote{It is also called 
the unit-scale randomization in \cite{DLM}.} of $u_0$ by
\begin{equation}
u_0^\omega  := \sum_{n \in \Z^d} g_n (\omega) \psi(D-n) u_0.
\label{gauss4}
\end{equation}

\noi
Compare this with the randomization \eqref{gauss3} on $\T^d$.
Under the assumption \eqref{R2}, Lemmas \ref{LEM:Hs}
and \ref{LEM:PStr1} also hold for the Wiener randomization \eqref{gauss4}
of a given function on $\R^d$.
The proofs remain essentially the same
with an extra ingredient of the following version of Bernstein's inequality:
\begin{equation}
\|\psi(D -n) u_0 \|_{L^q(\R^d)} \les \|\psi(D-n)  u_0  \|_{L^p(\R^d)}, \qquad
1\leq p \leq q \leq \infty,
\label{R4}
\end{equation}

\noi
for all $n \in \Z^d$.
The point of \eqref{R4} is that once a function is (roughly) restricted to a unit cube, 
we do not need to lose any derivative to go from the $L^q$-norm  to the $L^p$-norm, $q \geq p$.
See \cite{BOP1} for the proofs
of the analogues of Lemmas \ref{LEM:Hs}
and \ref{LEM:PStr1}.

Note that the probabilistic Strichartz estimates in Lemma \ref{LEM:PStr1}
hold only locally in time.
On $\T^d$, this does not cause any loss since the usual deterministic Strichartz estimates
also hold only locally in time.
On the other hand, 
the  Strichartz estimates on $\R^d$ hold globally in time:
\begin{equation}
\| e^{-it \Dl} u_0 \|_{L^q_t L^r_x (\R\times \R^d)} \lesssim \|u_0\|_{L^2_x(\R^d)}
\label{Str1}
\end{equation}

\noi
for any Schr\"odinger admissible pair
$(q, r)$, satisfying 
\begin{equation}
\frac{2}{q} + \frac{d}{r} = \frac{d}{2}
\label{Str2}
\end{equation}

\noi
 with $2\leq q, r \leq \infty$
and $(q, r, d) \ne (2, \infty, 2)$.
By incorporating the global-in-time estimate \eqref{Str1}, 
one can also obtain the following global-in-time probabilistic Strichartz estimates.

\begin{lemma}
Given $u_0 \in L^2(\R^d)$,
let $u_0^\o$ be its Wiener randomization defined in \eqref{gauss4}, satisfying \eqref{R2}.
Given a Schr\"odinger admissible pair $(q, r)$ with $q, r < \infty$,
let $\wt {r} \geq r$.
Then, there exist $C, c>0$ such that
\begin{align*}
P\Big( \|e^{-it\Dl} u_0^\omega\|_{L^q_t L^{\wt{r}}_x ( \R \times \R^d)} > \ld\Big)
\leq Ce^{-c \ld^2 \|u_0\|_{L^2}^{-2}}
\end{align*}

\noi
for all $\ld > 0$.

\end{lemma}

As in the periodic setting, 
similar global-in-time probabilistic Strichartz estimates
also hold
for the wave equation.
See  \cite{LM,  POC, OP}.

\begin{remark}\rm
(i) As we point out below, the Wiener randomization is special among other possible randomizations stemming from functions spaces in time-frequency analysis.
Recall the following definition of the modulation spaces in time-frequency analysis
\cite{Fei, FG1, FG2}.
Let $0<p,q\leq \infty$ and $s\in\mathbb R$; $M^{p, q}_s$ consists of all tempered distributions $u\in\mathcal S'(\R^d)$ for which the (quasi) norm
\begin{equation*}
\|u\|_{M_s^{p, q}(\R^d)} := \big\| \jb{n}^s \|\psi(D-n) u
\|_{L_x^p(\R^d)} \big\|_{\l^q_n(\mathbb{Z}^d)}
%\label{mod2}
\end{equation*}

\noi
is finite, where $\psi(D-n)$ is as above.
In particular, we see that 
the Wiener randomization \eqref{gauss4}
based on the unit cube decomposition of the frequency space
is 
very natural from the perspective of time-frequency analysis associated with the modulation spaces.

 \smallskip
 
 \noi
 (ii)
Let $\varphi_0, \varphi \in \mathcal{S}(\R^d)$ such that
$\supp \varphi_0 \subset \{ |\xi| \leq 2\}$, $\supp \varphi \subset
\{ \frac{1}{2}\leq |\xi| \leq 2\}$, and $ \varphi_0(\xi) + \sum_{j =
1}^\infty \varphi(2^{-j}\xi) \equiv 1.$ 
With $\varphi_j(\xi) =\varphi(2^{-j}\xi)$,
one may consider the following decomposition of a function:
\begin{equation}
u_0 = \sum_{j = 0}^\infty  \varphi_j(D) u_0
\label{gauss5}
\end{equation}

\noi
and introduce
the following randomization:
\begin{equation*}
u_0^\omega : = \sum_{j = 0}^\infty g_n (\omega) \varphi_j(D) u_0.
\end{equation*}

\noi
Note that \eqref{gauss5} can be viewed
as a  decomposition 
associated with the Besov spaces.
In view of the Littlewood-Paley theory,
such a randomization does not yield any improvement
on differentiability or integrability
and thus it is of no interest.

\smallskip
 
 \noi
 (iii)
Consider 
the following  wavelet series of a function:
\begin{equation}
u_0 = \sum_{\ld \in \Ld}  \jb{u_0, \psi_\ld} \psi_\ld, 
\label{gauss6}
\end{equation}

\noi
where $\{\psi_\ld\}_{\ld\in \Ld}$ is a wavelet basis of $L^2(\R^d)$.
One may also fancy  the following  randomization based on 
the   wavelet expansion \eqref{gauss6}:
\begin{equation}
u_0^\o := \sum_{\ld \in \Ld} g_\ld(\o) \jb{u_0, \psi_\ld} \psi_\ld.
\label{gauss7}
\end{equation}

\noi
Under some regularity assumption on $\psi_\ld$, 
we have the following characterization of the $L^p$-norm \cite[Chapter 6]{Meyer}:
\begin{align}
\| u_0^\o\|_{L^p(\R^d)} 
\sim \bigg\|\bigg( \sum_{\ld \in \Ld}
|g_\ld(\o)|^2 |\jb{u_0, \psi_\ld}|^2 |\psi_\ld(x)|^2
\bigg)^\frac{1}{2}\bigg\|_{L^p(\R^d)}
\label{gauss8}
\end{align}

\noi
for $1 < p < \infty$.
For example, if $\{ g_\ld\}_{\ld \in \Ld}$ is a sequence of independent Bernoulli random variables, 
then it follows from \eqref{gauss8} that
$\|u_0^\o\|_{L^p} \sim \|u_0\|_{L^p}$
and hence we see no improvement on integrability
under the randomization \eqref{gauss7}.

\end{remark}

\section{Probabilistic well-posedness of NLW and NLS}
\label{SEC:WP}

In this section, we go over some aspects of probabilistic well-posedness of nonlinear dispersive PDEs.
In recent years, there has been an increasing number of 
 probabilistic well-posedness results for these equations. 
See \cite{BOP2, POC, OTz, OOP} 
and the references therein.

\subsection{Basic almost sure local well-posedness argument}\label{SUBSEC:LWP}
In the following, we consider the defocusing\footnote{For the local-in-time argument, 
the defocusing/focusing  nature of the equation does not play any role.} cubic NLW on $\T^3$:
\begin{equation}
\begin{cases}\label{NLW2}
 \dt^2  u = \Delta u - u^3  \\
(u, \dt u) |_{t = 0} = (u_0, u_1) \in \H^s(\T^3), 
\end{cases}
\qquad ( t, x) \in \R \times \T^3.
\end{equation}

\noi
We say that $u$ is a solution to \eqref{NLW2} if $u$ satisfies the following Duhamel formulation:
\begin{equation}
u(t)=S(t)(u_0, u_1)
-\int_{0}^t \frac{\sin ((t-t')|\nabla|)}{|\nabla|}u^3(t')dt', 
\label{NLW3}
\end{equation}

\noi
where $S(t)$ denotes the linear wave operator:
\begin{equation*}
S(t)\left(u_0, u_1\right):=\cos(t|\nabla|)u_0+\frac{\sin (t|\nabla|)}{|\nabla|}u_1.
\end{equation*}

\noi
In view of  \eqref{scaling2}, 
we see that the scaling-critical regularity for \eqref{NLW2} is $s_\text{crit} = \frac 12$.
When $s\geq \frac 12$, 
 local well-posedness of \eqref{NLW2} in $\H^s(\T^3)$
 follows from a standard fixed point argument with the 
 (deterministic) Strichartz estimates.
Moreover, the equation \eqref{NLW2} is known to be ill-posed
in $\H^s(\T^3)$ for $s < \frac 12$ \cite{CCT, BT1, Xia}.
In the following, we take initial data $(u_0, u_1)$ to be in 
$ \H^s(\T^3) \setminus \H^{\frac 12}(\T^3)$
for some appropriate $s < \frac 12$.

Fix  $(u_0, u_1)\in \H^s(\T^3) \setminus \H^{\frac 12}(\T^3)$.
We apply the randomization defined in \eqref{gauss3}
to $(u_0, u_1)$.
More precisely, 
we set 
\begin{align}
(u_0^\o, u_1^\o)(x) 
:=  \bigg( \sum_{n \in \Z^3} g_{0,  n}(\o)\ft u_0(n)e^{in \cdot x}, 
 \sum_{n \in \Z^3} g_{1,  n}(\o)\ft u_1(n)e^{in \cdot x}\bigg), 
\label{gauss9}
 \end{align}

\noi
where 
$\{ g_{j, n} \}_{j = 0, 1, n \in \Z^3}$ is a sequence of independent mean-zero 
complex-valued random variables
conditioned that $g_{j, -n} = \cj{g_{j, n}}$, $j = 0, 1$, $n \in \Z^3$.
Moreover, we assume the exponential moment bound of type \eqref{R2}.

\begin{theorem}\label{THM:LWP1}
Let $s \geq 0$. Then, the cubic NLW \eqref{NLW2} on $\T^3$ is 
almost surely locally well-posed with respect to the randomization \eqref{gauss9}.
Moreover, the solution 
 $u$ lies in the class:
\begin{align*}
S(t)(u_0^\o, u_1^\o) + C([0, T_\o]; H^1(\T^3))
\subset C([0,T_\o] ; L^2(\T^3))
\end{align*}

\noi
for $T_\o  = T_\o(u_0, u_1)>0$ almost surely.

\end{theorem}
	
This theorem is implicitly included in \cite{BT3}.
See also \cite{POC}.
We point out, however, that the main goal of the paper \cite{BT3} 
is to establish almost sure global well-posedness (see the next subsection)
and to introduce the notion of  probabilistic continuous dependence.
See \cite{BT3} for details.

In view of the Duhamel formulation \eqref{NLW3}, 
we write the solution $u$ as 
\begin{align}
u = z + v,
\label{NLW4a}
\end{align}

\noi
where $z = z^\o = S(t) (u_0^\o, u_1^\o)$
denotes the random linear solution.
Then, instead of studying the original equation \eqref{NLW2}, we study the equation satisfied by 
the nonlinear part $v$:
\begin{equation}
\begin{cases}\label{NLW5}
 \dt^2  v = \Delta v - (v+z^\o)^3  \\
(v, \dt v) |_{t = 0} = (0, 0)
\end{cases}
\end{equation}

\noi
by viewing the random linear solution $z^\o$ as an explicit external source term.
Given $\o \in \O$, define $\G^\o$ by 
\begin{equation}
\G^\o (v) (t)=
-\int_{0}^t \frac{\sin ((t-t')|\nabla|)}{|\nabla|}(v+z^\o)^3(t')dt'.
\label{NLW6}
\end{equation}

\noi
Our main goal is to show that 
\begin{align}
v = \G^\o(v)
\label{NLW6a}
\end{align}
on some random time interval $[0, T_\o]$ with $T_\o > 0$ almost surely.\footnote{Needless
to say, the solution $v$ is random
since it depends on the random linear solution $z^\o$.  For simplicity, however, we suppress the superscript $\o$.}
Then, the solution $u$ to \eqref{NLW2} with the randomized initial data
$(u_0^\o, u_1^\o)$ in \eqref{gauss9} is given by \eqref{NLW4a}.

Given $T > 0$, we use the following shorthand notations: $C_T B = C([0, T]; B)$
and $L^q_TB = L^q([0, T]; B)$.
We also denote by  $B_1 \subset C_T \dot H^1$ the unit ball
in $C_T \dot H^1$ centered at the origin.
Suppose that $(u_0, u_1) \in \H^0(\T^3)$.
Then, by Sobolev's inequality, we have
\begin{align*}
\| \G^\o (v) \|_{C_T\dot H^1}
\leq \|v+z^\o\|_{L^3_T L^6_x}^3
& \le C_1  
 T \|v\|_{C_T\dot H^1}^3 +C_2 \|z^\o\|_{L^3_T L^6_x}^3\notag\\
& \le C_1  
 T \|v\|_{C_T\dot H^1}^3 + \tfrac 12,
\end{align*}

\noi
where the last inequality holds 
on  a set $\O_T $
thanks to the probabilistic Strichartz estimate
(i.e.~an analogue to Lemma \ref{LEM:PStr1} for the linear wave equation).
Moreover, we have 
\begin{align}
P(\O_T^c) \to 0  \text{ as }T \to 0.
\label{NLW7a} 
\end{align}

\noi
By taking $T> 0$ sufficiently small, we conclude that 
\begin{align}
\| \G^\o (v) \|_{C_T\dot H^1} \leq 1
\label{NLW7b}
\end{align}

\noi
for any $v \in B_1$ and $\o \in \O_T$.

Similarly, by taking $T > 0$ sufficiently small,  the following difference estimate holds:
\begin{align}
\|\G^\o (v_1) -  \G^\o (v_2)  \|_{C_T\dot H^1}%\notag\\
&  \le
C T^\frac{1}{3}\Big(  \sum_{j = 1}^2 T^\frac{2}{3}\|v_j\|_{C_T\dot H^1}^2
 + \|z^\o\|_{L^3_T L^6_x}^2\Big) 
 \|v_1 - v_2\|_{C_T\dot H^1}\notag\\
&  \le
\tfrac 12 
 \|v_1 - v_2\|_{C_T\dot H^1}
\label{NLW8}
\end{align}

\noi
for any $v_1, v_2 \in B_1$ and $\o \in \O_T$.
Therefore, from \eqref{NLW7b} and \eqref{NLW8}, 
we conclude that, given $T > 0$ sufficiently small, 
$\G^\o$ is a contraction on $B_1$ for any $ \o \in \O_T$.
By the fundamental theorem of calculus: $v(t) = \int_0^t \dt v(t') dt'$
and Minkowski's inequality, 
we also conclude that $v \in C([0, T];  H^1(\T^3))$.

Now, set $\Si = \bigcup_{0< T\ll1} \O_T$. 
Then, for each $\o \in \Si$, there exists a unique solution $v$
to \eqref{NLW6a} on $[0, T_\o]$ with $T_\o > 0$.
Moreover, it follows from \eqref{NLW7a} that $P(\Si) = 1$.
This proves Theorem \ref{THM:LWP1}.

The main point of the argument above
is (i) the decomposition \eqref{NLW4a}
and (ii) the gain of integrability for the random linear solution $z$ thanks to the probabilistic Strichartz estimates.
Then, the gain of one derivative in the Duhamel integral operator in \eqref{NLW6}
allows us to conclude the desired result in a straightforward manner.
Lastly, note that the same almost sure local well-posedness result holds
for the cubic NLW posed on $\R^3$ with the verbatim proof.

\begin{remark}\rm
 The decomposition~\eqref{NLW4a} allows us to separate
 the unknown part~$v$ from the explicitly known random part $z$
 and exploit the gain of integrability on $z$.
This  idea goes back to the work of 
McKean \cite{McKean} and 
Bourgain \cite{BO96}. 
See also Burq-Tzvetkov \cite{BT1}.
In the field of stochastic parabolic  PDEs, this decomposition
is usually referred to as the Da Prato-Debussche trick \cite{DPD}.

\end{remark}

Next, let us briefly discuss the situation for the cubic NLS on $\R^3$:
\begin{equation}
\begin{cases}
i \partial_t u = \Delta u - |u|^2 u  \\
u|_{t = 0} = u_0 \in H^s(\R^3),
\end{cases}
\qquad ( t, x) \in \R \times \R^3.
\label{NLS2}
\end{equation}

\noi
As before, 
the scaling-critical regularity for \eqref{NLS2} is $s_\text{crit} = \frac 12$
and  local well-posedness of \eqref{NLS2} in $H^s(\R^3)$, $s\geq \frac 12$, 
 follows from a standard fixed point argument with the Strichartz estimates.
Moreover, the equation \eqref{NLS2} is known to be ill-posed
in $H^s(\R^3)$ for $s < \frac 12$.
Nonetheless, we have the following almost sure local well-posedness.

\begin{theorem} \label{THM:LWP2}
Let $\frac 14 < s < \frac{1}{2}$.
Given $u_0 \in H^s (\R ^3)$, let $u_0^{\omega}$ 
be its Wiener randomization defined in \eqref{gauss4}, satisfying \eqref{R2}.
Then, 
the cubic NLS \eqref{NLS2} on $\R^3$ is almost surely locally well-posed
with respect to the random initial data $u_0^\o$.
Moreover, 
the solution $u$ lies in the class:
\begin{align*}
 e^{-it\Dl} u_0 ^{\omega} + C([0,T] ; H^\frac{1}{2}(\R^3)) \subset C([0,T] ; H^s(\R^3))
\end{align*}

\noi
for $T_\o  = T_\o(u_0)>0$ almost surely.
\end{theorem}

As in the case of the cubic NLW, write the solution $u$ as 
\begin{align}
u = z + v,
\label{NLS2a} 
\end{align}

\noi
where $z$ denotes the random linear solution:
\begin{align}
z = z^\o = e^{-it \Dl}u_0 ^{\omega}.
\label{NLS3}
\end{align}

\noi
Then, the residual part $v$ satisfies the following perturbed NLS:
\begin{equation}
\begin{cases}
	 i \dt v  =  \Dl v -  |v + z|^2(v+z)\\
v|_{t = 0} = 0.
 \end{cases}
\label{NLS4}
\end{equation}

\noi
In terms of the Duhamel formulation, we have
\begin{equation}
v(t) =  i \int_0^t e^{-i (t-t')\Dl}	  |v + z|^2(v+z)(t') dt'.
\label{NLS5}
\end{equation}

\noi
A key difference from  \eqref{NLW6}
is that there is no explicit smoothing on the Duhamel integral operator
in \eqref{NLS5}.
Hence, we need another mechanism to gain derivatives.
For this purpose, we employ the Fourier restriction norm method introduced by Bourgain \cite{BO1}.
The basic strategy for proving Theorem \ref{THM:LWP2}
is to expand the product $ |v + z|^2(v+z)$ in \eqref{NLS5}
and carry  out  case-by-case analysis
on terms of the form: 
\begin{align}
 v \cj v v, \ \ v \cj v z, \ \ v \cj z z,\ \  \dots, \ \ z\cj z z.
\label{NLS5a}
 \end{align}

In the following, we describe the main idea of the argument.
See \cite{BOP2} for the full details.
By the duality, it suffices to estimate
\begin{align}
 \bigg| \int _0^T \int_{\R^3} \jb{\nb}^\frac{1}{2}(v_1 \cj v_2 v_3) \cj v_4 dx dt\bigg|,
 \label{NLS6}
\end{align}

\noi
where $v_j = v$ or $z$, $j = 1, 2, 3$, and $v_4$ denotes the duality variable
at the spatial regularity 0.
By applying the dyadic decompositions to each function in \eqref{NLS6}, 
we separate the argument into several cases.
For the sake of the argument, let us denote by $N_j$  the dyadic size of the spatial frequencies
of $v_j$ after the dyadic decomposition
and assume $N_1 \sim N_4 \geq N_2 \geq N_3$.\footnote{Note that, due to the spatial integration in \eqref{NLS6}, 
the largest two frequencies of the dyadic pieces must be comparable.}
(i) When $N_2, N_3 \ges N_1^\theta$ for some suitable $\theta = \theta(s)\in (0, 1)$, 
we can move the derivatives from $v_1$ to  $v_2$ and $v_3$.
(ii) When $N_2, N_3 \ll N_1^\theta$, 
we group $v_1 \cj v_2$ and $v_3\cj v_4$ and apply 
the bilinear refinement of the Strichartz estimate 
\cite{BO98, OzawaT}, which allows us to gain some derivative.
(iii)~When $N_2\ges  N_1^\theta \geq N_3$,  
we combine (i) and (ii).
This allows us to prove Theorem \ref{THM:LWP2}.

\begin{remark}\rm
(i) Note that we did not need to perform 
any refined case-by-case analysis in proving Theorem \ref{THM:LWP1}
for the cubic NLW \eqref{NLW2}.
This is thanks to the explicit gain of one derivative in the Duhamel integral operator \eqref{NLW6}.

\smallskip

\noi
(ii) One may also study almost sure local well-posedness of the cubic NLS
posed on $\T^3$ below the scaling-critical regularity.
In this case, the argument becomes more involved due to the lack of 
the bilinear refinement of the Strichartz estimate,
which was the main tool for gaining derivatives in the problem on $\R^3$.
See for example \cite{BO96, CO, NS}. 
We point out that, in these works, 
the random initial data was taken to be of the specific form \eqref{gauss1}.
Then, the main ingredient is 
the probabilistic improvement of Young's inequality (Lemma \ref{LEM:PY}), 
which allows us to save some summations.
Note that with the random initial data of the form \eqref{gauss1}, 
this probabilistic improvement on summability allows us to gain derivatives
thanks to the reciprocal powers of the spatial frequencies in \eqref{gauss1}.
At this point, however, 
for the general randomized initial data \eqref{gauss3} on $\T^d$, 
it is not clear to us how to prove almost sure local well-posedness below the critical regularity.

\end{remark}

\begin{remark}\rm
(i) The solution map: $(u_0, u_1) \mapsto u$ for \eqref{NLW2} is   classically ill-posed  when  
$s < s_\text{crit} = \frac  12$.
The decomposition \eqref{NLW4a} tells us that we can decompose the ill-posed solution map 
as 
\begin{align}
(u_0^\o, u_1^\o) \stackrel{\footnotesize \circled{1}}{\longmapsto} z^\o 
\stackrel{\footnotesize \circled{2}}{\longmapsto} v
\longmapsto u = z^\o + v, 
\label{decomp1}
\end{align}

\noi
where the first step 
$\small \circled{1}$
involves stochastic analysis (i.e.~the probabilistic Strichartz estimates)
and the second step 
$\small \circled{2}$ is entirely deterministic. 
Moreover, there is continuous dependence of the nonlinear part $v$
on the random linear part $z$ in the second step.
See Remark \ref{REM:BB} for more on this issue.

\smallskip
\noi
(ii) By a similar argument, one can prove almost sure local well-posedness
with  initial data of the form:~``a smooth deterministic function 
$+$ a rough random perturbation''.
For example,  
 given deterministic $(v_0, v_1) \in \H^1(\T^3)$, 
consider the random initial data
$(v_0, v_1) + 
(u_0^\o, u_1^\o)$
for \eqref{NLW2},  where  $(u_0^\o, u_1^\o)$ is as in \eqref{gauss9}.
Then, the only modification in the proof of Theorem~\ref{THM:LWP1} appears %argument above appears 
in that the initial data for the perturbed equation \eqref{NLW5}
is now given by $(v_0, v_1)$.
In this case, we have the following decomposition of the ill-posed solution map:
\begin{align*}
(v_0 + u_0^\o, v_1 + u_1^\o) \longmapsto (v_0, v_1, z^\o) 
\longmapsto v
\longmapsto u = z^\o + v.
\end{align*}

\noi
See \cite{OQ} for a further discussion on the random initial data
of this type.

\smallskip
\noi
(iii) In Theorems \ref{THM:LWP1} and \ref{THM:LWP2}, 
we used the term ``almost sure local well-posedness'' in a loose manner,
following Bourgain.
In fact, 
what is claimed  in Theorems \ref{THM:LWP1} and \ref{THM:LWP2}
is simply almost sure local existence of unique solutions.
We point out, however, that the decompositions \eqref{NLW4a}
and \eqref{NLS2a} provide
continuous dependence of the nonlinear part $v$ (in a higher regularity)
on the random linear solution $z$ as mentioned above.

\smallskip

\noi
(iv) 
In \cite{BT3}, Burq-Tzvetkov 
introduced the notion of probabilistic continuous dependence
of $u$ on the random initial data, 
thus providing a more complete notion of probabilistic well-posedness.

\smallskip

\noi
(v) Unlike the usual deterministic theory, 
the approximation property of the random solution $u^\o$ constructed
in Theorem \ref{THM:LWP1} by smooth solutions
crucially depends on a method of approximation.
On the one hand, Xia \cite{Xia} showed that
the solution map: $(u_0, u_1) \mapsto u$ for \eqref{NLW2} 
is discontinuous everywhere in $\H^s(\T^3)$ when $s < \frac{1}{2}$.
This in particular shows that the solution map for \eqref{NLW2}, a priori  defined on smooth functions, 
does not extend continuously to rough functions, including the case of the random initial data
$(u_0^\o, u_1^\o)$ considered in 
Theorem \ref{THM:LWP1}.\footnote{In fact, almost sure norm inflation
at $(u_0^\o, u_1^\o)$ holds.}
On the other hand, 
by considering mollified initial data  $(\eta_\eps * u_0^\o, \eta_\eps * u_1^\o)$, 
one can show that the corresponding smooth solution $u_\eps$
converges almost surely to the solution $u^\o$ constructed in Theorem \ref{THM:LWP1}.
Moreover, the limit is independent of the mollification kernel~$\eta$.
See \cite{TzNote, Xia, OOTz}.

\smallskip

\noi
(vi) 
In \cite{OOP}, Oh-Okamoto-Pocovnicu proved almost sure local well-posedness
of  the following NLS without gauge invariance:
\begin{equation*}
i \partial_t u = \Delta u - |u|^p
\end{equation*}

\noi
in the regime where non-existence of solutions is known.
This shows an example of a probabilistic argument overcoming a stronger form of ill-posedness
than discontinuity of a solution map (which is the case for NLW \eqref{NLW2} and NLS \eqref{NLS2}
discussed above).
In the same paper, the authors also discussed a probabilistic construction of finite time blowup solutions
below the scaling-critical regularity.

\smallskip

\noi
(vii) In this section, we discussed almost sure local well-posedness 
based on a simple Banach fixed point argument.
There are also 
 probabilistic constructions of local-in-time solutions
 which are not based on a contraction argument. 
 See \cite{OhFE, Richards, OTW}.
 See also 
Subsection \ref{SUBSEC:BB}.

\end{remark}

\subsection{On almost sure global well-posedness}
\label{SUBSEC:GWP}

Before going over some of the almost sure global well-posedness results
in the literature, we point out that, in the stochastic setting, it suffices to prove
the  following statement to conclude  almost sure global well-posedness.

\begin{lemma}
[``Almost'' almost sure global well-posedness]
\label{LEM:GWP}
Given  $T > 0$ and  $\eps > 0$, there exists $ \Omega_{T, \eps}\subset \O$
with  $P(\O_{T, \eps}^c)< \eps$ 
such that a solution $u^\o$ exists on $[-T, T]$
for any $\o \in \O_{T, \eps}$.

\end{lemma}

\noi
Then, almost sure global well-posedness follows from  Lemma \ref{LEM:GWP} and Borel-Cantelli lemma.

In \cite{BO94}, Bourgain studied the invariance property of  the Gibbs 
measure $\rho$ in \eqref{Gibbs2}
for NLS \eqref{NLS1} on $\T$.
The main difficulty of this problem is the construction of 
global-in-time dynamics in the support of the Gibbs measure, 
i.e.~in $H^\s(\T)$, $\s < \frac 12$.
By introducing the Fourier restriction norm method, 
Bourgain \cite{BO1} proved local well-posedness of NLS 
below $H^\frac{1}{2}(\T)$.
Global well-posedness, however, was obtained only for the cubic case ($p = 3$).
In~\cite{BO94}, Bourgain combined PDE analysis with ideas from probability and dynamical systems
and proved global well-posedness of NLS almost surely with respect to the Gibbs measure~$\rho$.
The main idea is to use the (formal) invariance of the Gibbs measure $\rho$
as a replacement of a conservation law, providing a control on the growth of the relevant norm
of solutions in a probabilistic manner.
More precisely, he exploited the invariance of the truncated Gibbs measure $\rho_N$
associated with the truncated NLS:
\begin{equation}
\label{NLS7}
i \partial_t u_N =  \Delta u_N - \P_N(|\P_N u_N|^{p-1} \P_N u_N),
\end{equation}

\noi
where $\P_N$ denotes the Dirichlet projection onto the 
frequencies $\{|n|\leq N\}$, and proved the following growth bound; 
given $N \in \NB$, $T > 0$,  and  $\eps > 0$, there exists $ \Omega_{N, T, \eps}\subset \O$
with  $P(\O_{N, T, \eps}^c)< \eps$ 
such that for $\o \in \O_{N, T, \eps}$, 
the solution $u_N^\o$ to \eqref{NLS7} satisfies
\begin{align}
 \| u_N(t) \|_{H^\s} \leq C \Big( \log \frac{T}{\eps} \Big)^\frac{1}{2}
\label{NLS8}
\end{align}

\noi
 for any $t \in [-T, T]$, where $C$ is independent of $N$.
Combining \eqref{NLS8} with a standard PDE analysis, 
the same estimate also holds for the solution $u$ to \eqref{NLS1}
on a set $\O_{T, \eps}$ with $P(\O_{T, \eps}^c) < \eps$, 
yielding Lemma \ref{LEM:GWP}.
We point out that invariance of the Gibbs measure $\rho$ follows easily once
we have well defined global-in-time dynamics.
This argument is now known as {\it Bourgain's invariant measure argument}
and is widely applied\footnote{This argument 
 is not limited to nonlinear dispersive PDEs.
For instance, see  \cite{HairerM} for an application of this argument in studying 
 a stochastic parabolic PDE.} in situations where there is a formally invariant measure.

Note that while Bourgain's invariant measure argument
is very useful, its use is restricted to the situation
 where there is a formally invariant measure.
Namely, 
it can not be used to 
study the global-in-time behavior of solutions to an evolution equation 
with general random initial data.
In the following, we list various methods in establishing
almost sure global well-posedness of nonlinear dispersive PDEs
with general random initial data.

\smallskip

\noi
$\bullet$
{\it Bourgain's high-low method in the probabilistic context}:
In \cite{BO98}, Bourgain introduced an argument
to prove global well-posedness of NLS below the energy space.
The main idea is to divide the dynamics into the low 
and high frequency parts, where the low frequency part
lies in $H^1$ and hence the energy conservation is available.
The main ingredient in this argument is the nonlinear smoothing property of the high frequency part.
By exhibiting nonlinear smoothing in a probabilistic manner,\footnote{See
 Theorems \ref{THM:LWP1} and \ref{THM:LWP2} for
such nonlinear smoothing in the probabilistic setting.}
Colliander-Oh \cite{CO} implemented this argument in the probabilistic setting
to prove almost sure global well-posedness of the (renormalized) cubic NLS
on $\T$ in negative Sobolev spaces.
See also \cite{PRT, LM}.

\smallskip

\noi
$\bullet$
{\it Probabilistic a priori energy bound}:
The most basic way to prove global well-posedness
is to iterate a local well-posedness argument.
This can be implemented in a situation, where
one has  (deterministic) local well-posedness in the subcritical sense\footnote{Namely, 
the local existence time depends only on the norm of initial data.}
and the relevant norm of a solution is controlled by a conservation law.
In the probabilistic setting, one can implement a similar idea.
Burq-Tzvetkov \cite{BT3} proved almost sure global well-posedness
of the defocusing cubic NLW on $\T^3$
by estimating the growth of the (non-conserved) energy:
\[ H(v) = \frac{1}{2} \int_{\T^3} |\nb v|^2 dx+ 
 \frac{1}{2} \int_{\T^3} (\dt v)^2 dx
 +  \frac{1}{4} \int_{\T^3} v^4 dx\]

\noi
of the nonlinear part $v = u - z$, solving the perturbed NLW \eqref{NLW5}.
The argument is based on Gronwall's inequality
along with Cauchy-Schwarz' inequality
and the probabilistic Strichartz estimates.
There is a slight loss of derivative in controlling the $L^\infty_x$-norm
of the random linear solution and thus this argument works only for $s > 0$.
In the endpoint case ($s = 0$), Burq-Tzvetkov adapted Yudovich's argument \cite{Yu}
to control the energy growth of the nonlinear part $v$.

We point out that 
this argument based on Cauchy-Schwarz' and Gronwall's inequalities 
works only for the cubic case. 
In \cite{OP}, Oh-Pocovnicu proved almost sure global well-posedness
of the energy-critical defocusing quintic NLW on $\R^3$
(see the next item).
The main new ingredient is a new probabilistic a priori energy bound
on the nonlinear part $v = u - z$.
See \cite{LM2, SX} for  
almost sure global well-posedness of the three-dimensional NLW
$3 < p < 5$, following the idea in \cite{OP}.

\smallskip

\noi
$\bullet$
{\it Almost sure global existence by a compactness argument}:
In \cite{BTT1}, 
Burq-Thomann-Tzvetkov studied the defocusing cubic NLW \eqref{NLW2}
on $\T^d$, $ d\geq 4$.
By establishing a probabilistic energy bound on the nonlinear part $v_N$ of the solution to the truncated NLW,
which is uniform in the truncation parameter $N \in \NB$,  
they established a compactness property of $\{v_N\}_{N \in \NB}$, 
which allowed them to prove almost sure global existence.
See Nahmod-Pavlovi\'c-Staffilani \cite{NPS} for a precursor of this argument in the context of the Navier-Stokes equations.

In \cite{BTT2}, 
Burq-Thomann-Tzvetkov 
adapted a different kind of compactness argument in fluids by Albeverio-Cruzeiro \cite{AC}
to the dispersive setting 
and  constructed almost sure global-in-time dynamics
 with respect to Gibbs measures $\rho$ for various equations.
By exploiting the invariance of the truncated Gibbs measure $\rho_N$
for the truncated dynamics, they showed a compactness property
of the measures $\{ \nu_N \}_{N \in \NB}$ on space-time functions, 
where $\nu_N = \rho_N\circ \Phi_N^{-1}$ denotes the pushforward
of the truncated Gibbs measure $\rho_N$ under the global-in-time solution map $\Phi_N$ of the truncated dynamics.
Skorokhod's theorem then allowed them to construct 
a function $u$ as  an almost sure limit of the solutions $u_N$ to the truncated dynamics
(distributed according to $\nu_N$), 
yielding almost sure global existence.
See \cite{OTh, DSC} for more on this method.
This construction of  global-in-time solutions is  closely related
to the notion of martingale solutions
in the field of stochastic PDEs.
See~\cite{DZ}.

Due to the use of compactness, 
the global-in-time solutions constructed above 
are not 
unique (except for the two-dimensional case in \cite{NPS}).
In the four-dimensional case, 
the solutions to the defocusing cubic NLW constructed in \cite{BTT1}
were shown to be unique in \cite{OP2}, 
relying on the result \cite{POC}, 
which we discuss in the next item.

\smallskip

\noi
$\bullet$
{\it Almost sure global well-posedness in the energy-critical case via
perturbation/stability results}:
In the deterministic setting, the energy conservation allows us to prove
global well-posedness in $H^1$ of energy-subcritical defocusing NLW and NLS.
In the energy-critical setting, however, the situation is more complicated
and the energy conservation is not enough.
Over the last several decades, substantial effort was made in understanding global-in-time behavior of
solutions to the energy-critical NLW and NLS.
See \cite{BOP2, POC, OOP} for the references therein.

 In an analogous manner, a probabilistic a priori energy bound discussed above
does not suffice to prove almost sure global well-posedness
of the energy-critical defocusing NLW and NLS.
In \cite{BOP2}, we introduced a new argument, 
using  perturbation/stability results for NLS to approximate the dynamics 
of the perturbed NLS \eqref{NLS4}
by unperturbed NLS dynamics on short time intervals, 
and proved conditional almost sure global well-posedness
of the energy-critical defocusing cubic NLS on $\R^4$, 
provided that the energy of the nonlinear part $v = u - z$
remains bounded almost surely for each finite time.
Subsequently, by establishing probabilistic a priori energy bounds, 
Pocovnicu \cite{POC} and Oh-Pocovnicu \cite{OP, OP2}
applied this argument and 
proved almost sure global well-posedness of 
 the energy-critical defocusing NLW
 on $\R^d$ and $\T^d$, $d = 3, 4, 5$.
More recently, 
Oh-Okamoto-Pocovnicu
\cite{OOP}
established probabilistic a priori energy bounds
for the energy-critical defocusing NLS on $\R^d$,  $d = 5, 6$
and proved almost sure global well-posedness
for these equations.

\smallskip

\noi
$\bullet$
{\it Almost sure scattering results via a double bootstrap argument:}
The almost sure global well-posedness results mentioned above
do not give us any information on the qualitative behavior of global-in-time solutions such as scattering.
This is due to the lack of global-in-time space-time bounds in the argument mentioned above.
In a recent paper, 
Dodson-L\"uhrmann-Mendelson \cite{DLM}
studied the energy-critical defocusing NLW on $\R^4$ in the radial setting
and proved almost sure scattering in this setting.
The argument is once again  based on  applying 
perturbation/stability results as in  \cite{BOP2, POC, OP}.
The main new ingredient is the Morawetz estimate for the perturbed NLW \eqref{NLW5}
satisfied by the nonlinear part $v = u - z $.
More precisely, they implemented a double bootstrap argument,
controlling the energy and the Morawetz quantity for the nonlinear part $v$
in an intertwining manner.
Here, the radial assumption plays a crucial role in applying the radial improvement of the Strichartz estimates.
We also point out that even if the original initial data is radial,
its randomization is no longer radial and some care must be taken 
in order to make use of the radial assumption.
In \cite{KMV}, Killip-Murphy-Vi\c{s}an proved
an analogous almost sure scattering result for
 the energy-critical defocusing NLS on $\R^4$ in the radial setting.\footnote{We also mention
 a recent work by Dodson-L\"uhrmann-Mendelson \cite{DLM2} that appeared after the completion of this note.
 The main new idea in \cite{DLM2} is to adapt
 the functional framework for the derivative NLS and Schr\"odinger maps
to study the perturbed NLS \eqref{NLS4}. }
It would be of interest to remove the radial assumption imposed 
in the aforementioned work \cite{DLM, KMV}.

We also mention the work \cite{BTT0, DS1} on the almost sure scattering results
based on (a variant of) Bourgain's invariant measure argument
combined with the equation-specific transforms.

\smallskip

In this subsection, we went over various globalization arguments.
We point out, however, that,
except for Bourgain's invariant measure argument
(which is restricted to a very particular setting), 
all the almost sure globalization arguments are based
on known deterministic globalization arguments.
Namely, there is {\it no} 
probabilistic  argument at this point that is not based on a known deterministic argument
for generating global-in-time solutions.

\subsection{Further discussions}\label{SUBSEC:BB}

We conclude this section with a  further discussion on  probabilistic construction
of solutions to nonlinear dispersive PDEs.

\smallskip

\noi
$\bullet$ {\it Higher order expansions}:
In Subsection \ref{SUBSEC:LWP}, we described the basic
strategy for proving almost sure local well-posedness of NLW and NLS.
In the following, we  describe the main idea on how to improve the regularity thresholds
stated in Theorems \ref{THM:LWP1} and \ref{THM:LWP2}.

We first consider  the cubic NLS on $\R^3$.
The almost sure local well-posedness
stated in  Theorem \ref{THM:LWP2}
follows from the case-by-case analysis in \eqref{NLS5a}
within the framework of the Fourier restriction norm method.
See \cite{BOP2} for the details.
By examining the case-by-case analysis in \cite{BOP2}, 
we see that 
the regularity restriction $s >  \frac 14$ in Theorem \ref{THM:LWP2}
comes from 
the cubic interaction 
of the random linear solution:
\begin{align}
  z_3(t)  := i  \int_0^t S(t - t') |z_1|^2 z_1(t') dt',
\label{Z3}
\end{align}

\noi
where $z_1 : = z^\o = e^{-it \Dl} u_0^\o$ defined in \eqref{NLS3}.
On the one hand, given $u_0 \in H^s(\R^3)$, $ 0 \leq s < 1$, 
we can prove that $z_3$ in \eqref{Z3} has spatial regularity $2s-\eps$ for any $\eps > 0$.
On the other hand, we need $2s-\eps > s_\text{crit} = \frac 12$ in order to close the argument.
This yields the regularity restriction $s > \frac 14$ stated in Theorem \ref{THM:LWP2}.

By noting that all the other interactions in \eqref{NLS5a} behave better than $z_1 \cj{z_1} z_1$, 
we introduce the following second order expansion:
\begin{align*}
u = z_1 + z_3 + v
\end{align*}

\noi
to remove the worst  interaction $z_1 \cj{z_1} z_1$.
In this case,   the residual term $v := u - z_1 - z_3$ satisfies the following equation:
\begin{equation*}
\begin{cases}
	 i \dt v =  \Dl v -   \N(v + z_1+z_3) + \N(z_1)\\
v|_{t = 0} = 0,
 \end{cases}
\end{equation*}

\noi
where $\N(u) = |u|^2 u$.
In terms of the Duhamel formulation, we have
\begin{align}
v(t) & =  i \int_0^t e^{-i (t-t')\Dl} \big\{\N(v + z_1 + z_3) - \N(z_1)\big\}(t') dt'
\label{NLS10}
\end{align}

\noi
Note that we have removed the worst interaction $z_1\cj{z_1} z_1$
appearing in the case-by-case analysis \eqref{NLS5a}.
There is, however, a price to pay;
we need to carry out the following  case-by-case analysis
\begin{align}
v_1\cj{v_2} v_3, 
\quad \text{for  $v_i = v,$ $z_1$, or $z_3$, $i = 1, 2, 3$,  but not all $v_i$ equal to $z_1$}, 
\label{case1}
\end{align}

\noi
containing more terms than the previous case-by-case analysis \eqref{NLS5a}.
Nonetheless, 
 by studying the fixed point problem \eqref{NLS10} for $v$, 
 we can lower the regularity threshold for almost sure local well-posedness.
Note that the solution $u$ thus constructed lies in the class:
\begin{align*}
 z_1 + z_3 
+ C([0, T]; H^{\frac 12} (\R^3)) 
\subset C([0, T];H^s(\R^3))
\end{align*}

\noi 
for some appropriate $s < \frac 12$.

By examining the case-by-case analysis \eqref{case1}, 
we see that the worst interaction appears
in 
\begin{align}
\text{$  z_{j_1} \cj{z_{j_2}}z_{j_3} $
\quad with  $(j_1, j_2, j_3) = (1, 1, 3)$ up to permutations, }
\label{NLS11}
\end{align}
 giving rise to the following third order term:
\begin{align}
  z_5 (t) := i  \sum_{\substack{j_1 + j_2 + j_3 = 5\\j_1, j_2, j_3 \in \{1, 3\}  }}\int_0^t S(t - t') 
  z_{j_1} \cj{z_{j_2}}z_{j_3} (t') dt'.
\label{Z5}
\end{align}

\noi
A natural
next step is to 
remove this non-desirable  interaction in \eqref{NLS11}
in the case-by-case analysis in \eqref{case1} by 
considering 
 the following  third order expansion:
\begin{align*}
u = z_1 + z_3 + z_5+ v.
%\label{v3}
\end{align*}

\noi
In this case, the residual term $v := u - z_1 - z_3- z_5$ satisfies the following equation:
\begin{equation*}
\begin{cases}
\displaystyle
 i \dt v =  \Dl v -  \N(v + z_1+z_3+z_5) 
+   \sum_{\substack{j_1 + j_2 + j_3 \in \{3, 5\} \\j_1, j_2, j_3 \in \{1, 3\}  }}
 z_{j_1} \cj{z_{j_2}}z_{j_3}\\
v|_{t = 0} = 0
 \end{cases}
\end{equation*}

\noi
and thus we need to carry out 
the following case-by-case analysis:
\begin{align*}
v_1\cj{v_2} v_3
\quad & \text{for  $v_i = v, z_1, z_3$, or $z_5$, $i = 1, 2, 3$, such that }  \\
& \text{it is not of the form $z_{j_1}\cj{z_{j_2}}z_{j_3} $ with $j_1 + j_2 + j_3 \in \{3, 5\}$}.
\end{align*}

\noi 
Note the increasing number of combinations.
While it is theoretically possible to repeat this procedure based on partial power series expansions, 
it seems virtually impossible to  keep track of all the terms 
as the order of the expansion grows.
In \cite{BOP3}, we introduced modified higher order expansions to avoid this combinatorial nightmare
and established  improved almost sure local well-posedness
based on the modified higher order expansions of arbitrary length.

Next, we briefly discuss the cubic NLW \eqref{NLW2} on $\T^3$.
Theorem \ref{THM:LWP1} already provides almost sure local well-posedness
for $s \geq 0$
and hence we now need to take $s < 0$.
In this case, the cubic product $z^3$ appearing in the perturbed NLW \eqref{NLW5}
does not make sense and we need to {\it renormalize} the nonlinearity.
We assume that 
 the randomized initial data $(u_0^\o, u_1^\o) $ is of the form \eqref{gauss9}
with  $\ft u_j(n) = \jb{n}^{-\al+j}$, $j = 0, 1$, 
and that 
 $\{ g_{j, n} \}_{j = 0, 1, n \in \Z^3}$ in \eqref{gauss9} 
is a sequence of independent standard complex-valued Gaussian random variables
conditioned that $g_{j, -n} = \cj{g_{j, n}}$, $j = 0, 1$, $n \in \Z^3$.
Comparing this with \eqref{gauss1}, 
we see that $(u_0^\o, u_1^\o) $ 
is distributed according to the Gaussian measure $\mu_\al \otimes \mu_{\al - 1}$
supported on $\H^s(\T^3)$ for $s < \al - \frac 32$.
In this case, we can apply the Wick renormalization
(see \cite{GKO, OTh2})
and obtain the following renormalized equation for $v = u- z$:
\begin{align}
 \dt^2  v 
 & = \Delta v \, -: \! (v+z)^3 \!:  \notag\\
 & = \Delta v \, - v^3 - 3 v^2 z -3 v  : \!z^2\!: -  : \! z^3 \!:,  
\label{NLW9}
\end{align}

\noi
where $: \! z^\l \!:$ is the  standard Wick power of $z$,
having the spatial regularity $\l(\al - \frac 32)-\eps$, 
 $\l = 2, 3$.
By studying the equation \eqref{NLW9}, 
it is easy to prove almost sure local well-posedness
for $\al > \frac 43$, corresponding to $s > -\frac 16$.
Note that the worst term in \eqref{NLW9}
is given by 
$: \! z^3 \!:$ with the spatial regularity $3\al - \frac 92- \eps$.

As in the case of the cubic NLS, 
we shall consider the second order expansion:
\begin{align}
u = z_1 + z_3 + v, 
\label{decomp2}
\end{align} 

\noi
where 
$z_1 = z$ and 
\begin{equation*}
z_3(t) := 
-\int_{0}^t \frac{\sin ((t-t')|\nabla|)}{|\nabla|}:\!z_1^3 (t') \!:dt'.
%\label{NLW10}
\end{equation*}

\noi
Then, the residual term $v = u - z_1 - z_3$ satisfies
\begin{align}
 \dt^2  v 
 & = \Delta v \, -: \! (v+z_1 + z_3 )^3 \!:  + : \! z_1^3 \!:\notag\\
 & = \Delta v \, - (v+z_3)^3 - 3 (v+z_3)^2 z_1 -3 (v+z_3)  : \!z_1^2\!:.
\label{NLW11}
\end{align}

\noi
Note that the worst term $: \!z_1^3\!: $ in \eqref{NLW9} is now eliminated.
In \eqref{NLW11}, the worst contribution  
is given by $3 (v+z_3)  : \!z_1^2\!:$ with regularity $2\al - 3-\eps$.\footnote{
Here, we assume that $z_3$ has positive regularity.
For example, we know that  $z_3$ has  spatial regularity at least $3\al - \frac 92 + 1-\eps $
and hence $\al > \frac 76$ suffices.}
By solving the fixed point problem \eqref{NLW11} for $v$, we can lower the regularity threshold.
In this formulation, the regularity restriction arises in making sense of the
product $z_3 \, \cdot  : \!z_1^2\!:$ as distributions of regularities $3\al - \frac 72-\eps$ and $2\al - 3-\eps$.

\begin{remark}\label{REM:BB} \rm
The second order formulation \eqref{decomp2} 
yields the following decomposition of the ill-posed solution map:
\begin{align*}
(u_0^\o, u_1^\o)  \longmapsto (z_1,  : \!z_1^2\!:, z_3)
\longmapsto v
\longmapsto u = z_1 + z_3 + v.
\end{align*}

\noi
Compare this with \eqref{decomp1}.
As before,  the first step 
involves stochastic analysis, while 
the second step  is entirely deterministic. 
One can establish  a further improvement
to the argument sketched above,
providing a meaning to 
$z_3 \, \cdot  : \!z_1^2\!:$
in a probabilistic manner:\footnote{Recall the following  paraproduct decomposition
of the product $fg$ of two functions $f$ and $g$:
\begin{align*}
fg 
& = f\pl g + f \pe g + f \pg g\notag \\
 :\!& =\sum_{j < k-1} \varphi_j(D) f \, \varphi_k(D) g
+ \sum_{|j - k|  \leq 1} \varphi_j(D) f\,  \varphi_k(D) g
+ \sum_{k < j-1} \varphi_j(D) f\,  \varphi_k(D) g.
\end{align*}

\noi
Since the paraproducts $z_3 \, \pl \!: \!z_1^2\!:$ and $z_3 \, \pg \! : \!z_1^2\!:$ always make sense as distributions,
it suffices to give a meaning to the resonant product $z_3 \, \pe  \!: \!z_1^2\!:$
in a probabilistic manner.
}
\begin{align*}
(u_0^\o, u_1^\o)  \longmapsto (z_1,  : \!z_1^2\!:, z_3, z_3 \, \pe  \!: \!z_1^2\!:)
\longmapsto v
\longmapsto u = z_1 + z_3 + v.
\end{align*}

\noi
See \cite{OPT} for details.
%We will address this issue in a forthcoming paper.

\end{remark}

\begin{remark}\rm

In a recent paper \cite{PW}, Pocovnicu-Wang 
provided a simple argument
for constructing unique solutions
to NLS with random initial data
by exploiting  the dispersive estimate.
In the context of the cubic NLS \eqref{NLS2} on $\R^3$, 
the random initial data can be taken to be only in $L^2(\R^3)$.
Compare this with Theorem \ref{THM:LWP2}.
Their construction, however, places a solution $u$ only in the class:
\begin{align*}
 e^{-it\Dl} u_0 ^{\omega} + C([0,T] ; L^4(\R^3)), 
\end{align*}

\noi
which does not embed in 
$ C([0,T] ; H^s(\R^3))$.
See also \cite{OPW} for a related result
in the context of the stochastic NLS on $\R^d$.

\end{remark}

\smallskip

\noi
$\bullet$
{\it Bourgain-Bulut's argument}:
In Subsection \ref{SUBSEC:GWP}, 
we described a globalization argument when there is a formally invariant measure 
(Bourgain's invariant measure argument). 
We point out that this globalization argument requires a separate
 local well-posedness argument
 (whether deterministic or probabilistic).
In \cite{BB1, BB2, BB3}, Bourgain-Bulut presented
a new argument, where they exploited formal invariance
already in the  construction of local-in-time solutions.
In the following, we briefly sketch  the essential idea in \cite{BB1, BB2, BB3}
by taking NLS \eqref{NLS1} on $\T$
as an example.

The main goal is to show that the solution $u_N$ to the truncated NLS \eqref{NLS7}
converges to some space-time distribution $u$, which turns out to satisfy the original NLS \eqref{NLS1}.
This is done by exploiting the invariance of the truncated Gibbs measure $\rho_N$:
\begin{align*}
d\rho_N 
& = Z^{-1} e^{- \frac{1}{p+1} \int |\P_N u|^{p+1} } d\mu_1, 
\end{align*}

\noi
for the truncated equation \eqref{NLS7}.
Let $T \leq 1$.
Then, by the probabilistic Strichartz estimates (Lemma \ref{LEM:PStr1}), we have
\begin{align*}
\mu_1\Big( u_0: \| e^{-it \Dl} u_0\|_{L^q_t W^{\s, r}_x([0, T]\times \T)}> \ld\Big)
\leq C\exp (-c \ld^2)
\end{align*}

\noi
for any $\s < \frac 12$, 
 finite  $q \geq 2$, and $2 \leq  r \leq \infty$.
Then, by using the mutual absolute continuity between $\mu_1$ and $\rho_N$
and exploiting the invariance of $\rho_N$, we can upgrade this estimate to 
\begin{align}
\rho_N \Big( u_0: \| u_N \|_{L^q_t W^{\s, r}_x([0, T]\times \T)}> \ld\Big)
\leq C\exp (-c \ld^{c'}), 
\label{Gibbs6}
\end{align}

\noi
where $u_N$ is the solution to the truncated NLS \eqref{NLS7}
with $u_N|_{t = 0} = u_0$.
We stress that the constants in \eqref{Gibbs6} are independent of $N \in \NB$.

Given $M \geq N \geq 1$, let $u_M$ and $u_N$ be the solutions to \eqref{NLS7} with 
the truncation size $M$ and $N$, respectively.
Then, on a time interval $I_j = [t_j, t_{j+1}] \subset [0, T]$, we have 
\begin{align}
u_M(t)- u_N(t)
& =e^{-i (t-t_j)\dx^2 }(u_M(t_j) - u_N(t_j) ) \notag\\
& \hphantom{X}
+ i \int_{t_j}^t e^{-i (t-t')\dx^2}	(\P_M- \P_N)  |u_M|^{p-1}u_M(t') dt'
\notag \\
& \hphantom{X}
+ i \int_{t_j}^t e^{-i (t-t')\dx^2}	\P_N  (|u_M|^{p-1}u_M- |u_N|^{p-1}u_N)(t') dt'.
\notag \\
& = \1 + \II + \III.
\label{NLS13}
\end{align}

\noi
Fix $s < \frac{1}{2}$ sufficiently close to $\frac 12$.
We estimate each term on the right-hand side of \eqref{NLS13} in the $X^{s, b}$-norm with $b = \frac 12+$.
See \cite{TAO} for the basic properties of the $X^{s, b}$-spaces.
The first term $\1$ is trivially bounded by 
$\|u_M(t_j) - u_N(t_j) \|_{H^s}$.
Noting that the second term $\II$ is supported on high frequencies $\{|n|> N\}$, 
we can show that it tends to 0 as $N \to \infty$.
Here, we used \eqref{Gibbs6} with $s < \s < \frac 12$,
giving a decay of  $N^{s-\s}$ for $\|\II\|_{X^{s, b}([t_j, t_{j+1}])}$.
As for the last term $\III$, by the fractional Leibniz rule, \eqref{Gibbs6}, and the periodic Strichartz estimate,
we can estimate it by $\| u_M -u_N\|_{X^{s, b}([t_j, t_{j+1}])}$.\footnote{Here, the implicit constant depends
on the choice of $\ld$ in \eqref{Gibbs6}, which needs to be chosen in terms of $N$.
See \cite{BB2} for more on this issue.}
This allows us to iterate the argument on intervals $I_j$ 
to cover the entire interval $[0, T]$ and show that $\{u_N\}_{N \in \NB}$
is Cauchy in $X^{s, b}([0, T])$.

The almost sure local well-posedness argument presented in Subsection \ref{SUBSEC:LWP}
exploited the gain of integrability only at the level of the random linear solution.
The main novelty of the argument presented above is the use of invariance
in constructing local-in-time solutions, which 
yields the gain of integrability for the truncated solutions $u_N$,
uniformly in $N \in \NB$.

\section{Remarks and comments}

\noi
(i) In \eqref{scaling2}, we introduced the scaling-critical Sobolev regularity $s_\text{crit}$.
Note that this regularity  is based on the (homogeneous) $L^2$-based Sobolev spaces
$\dot H^s = \dot W^{s, 2}$.
Given $ 1\leq r \leq \infty$, we can also consider  
the scaling-critical Sobolev regularity
adapted to the $L^r$-based Sobolev spaces
$\dot W^{s, r}$:
\begin{align*}
 s_\text{crit}(r)  = \frac{d}{r} - \frac{2}{p-1}
\end{align*}

\noi
 such that the homogeneous $\dot{W}^{s_\textup{crit}(r), r}$-norm is invariant
under the dilation symmetry \eqref{scaling1}.
Heuristically speaking, the gain of integrability depicted in Lemmas \ref{LEM:Hs}
and \ref{LEM:PStr1} allows us to lower the critical regularity from 
$s_\text{crit} = s_\text{crit}(2) = \frac{d}{2} - \frac{2}{p-1}$
to $ s_\text{crit}(r)  =  - \frac{2}{p-1} +\eps$ for $r \gg 1$.\footnote{Things are not as simple as stated here due to 
the unboundedness of the linear solution operator on $L^r$, $r\ne 2$, 
for dispersive equations.
In the case of the nonlinear heat equation, however, this heuristics can be seen more clearly.
Consider the following nonlinear  heat equation on $\R^d$:
\begin{align}
 \dt u = \Dl u - |u|^{p-1}u
 \label{heat}
\end{align}

\noi
 with initial data $u_0 \in L^2(\R^d)$.
 In general, (when $4 < d(p-1)$ for example), we do not know how to construct a solution
 with initial data in $L^2(\R^d)$.
 By randomizing the initial data $u_0$ as in \eqref{gauss4}, 
 we see that the randomized initial data $u_0^\o$ lies almost surely in $L^r(\R^d)$
 for any finite $r\geq 2$.   Then, by taking $r > \frac{d(p-1)}{2}$, 
 we can apply the deterministic subcritical local well-posedness
 result in \cite{BC} to conclude (rather trivial) almost sure local well-posedness of \eqref{heat}
 with respect to the Wiener randomization $u_0^\o$.
This is an instance of ``making the problem subcritical'' by randomization.
  }
For example, 
in the cubic case ($p = 3$), we have  
$s_\text{crit}(r) = \frac{d}{r} - 1 \to -1$ as $ r \to \infty$,
which makes the problem considered in Subsection \ref{SUBSEC:LWP}
{\it subcritical} in some appropriate sense.

\smallskip

\noi
(ii) In recent years, there has been a significant development in the well-posedness theory of singular stochastic parabolic PDEs.
For example, 
the theory of regularity structures by Hairer \cite{Hairer} 
and the paracontrolled distributions introduced by Gubinelli-Imkeller-Perkowski \cite{GIP, CC}
allow us to make sense of  the following stochastic quantization equation (SQE, dynamical $\Phi^4_3$ model) on $\T^3$:
\begin{equation}
 \partial_t u =  \Delta u - u^3  + \infty\cdot u +  \xi, 
\label{SQE}
\end{equation}

\noi
where $\xi$ denotes the space-time white noise.
See also \cite{Kup}.
Moreover, it has been shown that the Gibbs measure $\rho$ in \eqref{Gibbs2}
(in a renormalized form)
is invariant under the dynamics of \eqref{SQE}
\cite{HairerM, AK}.
It would be of great interest to study a similar problem
for the defocusing cubic NLS and NLW on $\T^3$
with the Gibbs measure $\rho$ as initial data.
While the NLS problem seems to be out of reach (see (iii) below),
one may approach the NLW problem by adapting the paracontrolled calculus\footnote{At this point, 
we do not know how to apply the theory of regularity structures to study dispersive PDEs, partly because we do not know how to lift the Duhamel integral operator for dispersive PDEs to regularity structures.}
to the wave case.

\smallskip

\noi
(iii)  There has also been some development in the solution theory
for singular stochastic nonlinear dispersive PDEs
 \cite{OhSKdV, GKO}.
The problem of particular importance is the following
(renormalized) stochastic cubic nonlinear Schr\"odinger equation (SNLS) on $\T$
with additive space-time white noise forcing:
\begin{equation}
\label{SNLS}
i \partial_t u = \Delta u - |u|^{2} u  + 2 \infty \cdot u +  \xi.
\end{equation}

\noi
Since SNLS \eqref{SNLS} on $\T^d$ scales like \eqref{SQE} on $\T^d$, 
one may be tempted to think that they are of equal difficulty.
This, however, is completely false;
while SQE on $\T^d$, $d = 1, 2, 3$ is subcritical, 
SNLS on the one-dimensional torus $\T$ is {\it critical} in the following sense.

The linear heat semigroup $e^{t\Dl}$ is bounded on $L^\infty$
and thus the scaling-critical regularity for the cubic heat equation is given 
by $s_\text{crit}(\infty) = -1$ for any dimension.
The space-time white noise under the Duhamel integral operator:\footnote{This is 
the so-called  stochastic convolution.}
$\int_0^t e^{(t-t') \Dl} \xi(dt')$
has (spatial) regularity $-\frac{d}{2} + 1-\eps$.
Comparing these two regularities, we see that SQE \eqref{SQE}
is critical when $d = 4$ and is subcritical when $d = 1, 2, 3$. 
We point out that both the theory of regularity structures and the theory of paracontrolled distributions
are subcritical theories and can not handle \eqref{SQE} when  $d = 4$.

Similarly, 
by recalling that the linear Schr\"odinger  group $e^{-i t\Dl}$ is bounded on $L^2$
and is unbounded on any $L^r$, $ r\ne 2$, 
it seems reasonable to use $r = 2$
to compute 
the scaling-critical regularity for the cubic NLS, 
thus giving 
 $s_\text{crit}(2) = \frac{d}{2}-1$.
On the other hand, 
the stochastic convolution $\int_0^t e^{-i (t-t') \Dl} \xi(dt')$
in this case  does not experience any smoothing and thus 
has (spatial) regularity $-\frac{d}{2} -\eps$.
Comparing these two regularities, we see that SNLS \eqref{SNLS}
is critical already when $d = 1$.
We point out that the (deterministic) cubic NLS on $\T$ with the spatial white noise
as initial data, i.e.~\eqref{gauss1}
with $s = 0$, basically has the same difficulty.
This criticality is also manifested in the fact that the higher order iterates
such as $z_3$ in \eqref{Z3}
and $z_5$ in \eqref{Z5} do not experience any smoothing.\footnote{For example, 
for the subcritical SQE on $\T^3$, the second order iterate (an analogue of $z_3$ in \eqref{Z3}) gains one derivative
as compared to the stochastic convolution.}
Lastly, we mention a recent work
\cite{FOW}, where the second author (with Forlano and Y.\,Wang)
established local well-posedness of \eqref{SNLS}
with a slightly smooth noise $\jb{\dx}^{-\eps} \xi$, $\eps > 0$.

\smallskip

\noi
(iv)  Various methods and ideas in the random data Cauchy theory for nonlinear dispersive PDEs
are applicable to study stochastic nonlinear dispersive PDEs
thanks to the gain of integrability on a stochastic forcing term.\footnote{at least on $\T^d$.
On $\R^d$, there is a limitation on the gain of integrability.  See \cite{DD, OPW}.}
For example, the almost sure local well-posedness result
for  the (renormalized) cubic NLS in $H^s(\T)$, $s > -\frac 13$, by Colliander-Oh \cite{CO}
essentially implies local well-posedness of SNLS \eqref{SNLS} 
with  a smoothed noise 
 $\jb{\dx}^{-\eps} \xi$, $\eps > \frac 16$.

\smallskip

\noi
(v)  Thanks to Bourgain's invariant measure argument, 
we now have a good understanding of how to build an invariant measure
of Gibbs type based on a conservation law.
Note, however, that these measures are supported
on rough functions (except for completely integrable equations)
and we do not know how to construct invariant measures supported on smooth functions. 
Bourgain \cite{BO97} wrote
``the most important challenge is perhaps the question if we may produce an invariant measure which is supported 
by smooth functions.''

\smallskip

\noi
(vi) So far, we discussed how to construct solutions in a probabilistic manner.
It would be of interest to develop
 a probabilistic argument to get more qualitative information of solutions.
For example, the random initial data \eqref{gauss1}
lies almost surely in $W^{\s, \infty}(\T^d)$, $\s < s - \frac d2$.
On the one hand, the deterministic well-posedness theory
propagates only the $H^\s$-regularity of solutions.
On the other hand, quasi-invariance of $\mu_s$
\cite{TzBBM, OTz, OTzX, OTz2} implies that 
the $W^{\s, \infty}$-regularity is also propagated in an almost sure manner.
An interesting problem may be to use probabilistic tools
to study the growth of high Sobolev norms of solutions.
For example, the argument in \cite{TzBBM, OTz} provides
a probabilistic proof of 
a polynomial upper bound.

\begin{acknowledgment}

\rm 
\'A.\,B. is partially supported by a grant from the Simons Foundation (No. 246024).
T.\,O.~was supported by the European Research Council (grant no.~637995 ``ProbDynDispEq'').
The authors would like to thank Justin Forlano for careful proofreading.

\end{acknowledgment}

\end{document}